\documentclass[12pt]{article}
\usepackage{amsfonts, amsmath}
\usepackage{amssymb}
\newtheorem{theo}{Theorem}[section]
\newtheorem{pro}[theo]{Proposition}
\newtheorem{lem}[theo]{Lemma}
\newtheorem{cor}[theo]{Corollary}

\newcommand{\ccc}{{\cal C}}
\newcommand{\calX}{{\cal X}}

\def\norm#1{\left \|\, #1 \, \right \|}

 \def\supp{\mathop{\rm supp}}

 \newcommand{\ep}{\varepsilon}

\begin{document}

\title{Riesz transform on manifolds  and  Poincar\'e inequalities}

\date{}

\author{Pascal Auscher\footnote{
Research partially supported  by the European Commission (IHP
Network ``Harmonic
Analysis and Related Problems'' 2002-2006, Contract HPRN-CT-2001-00273-HARP). This research began at  the Centro De Giorgi of the Scuola Normale Superiore de Pisa on the occasion of a special program  of the network. The authors thank the organizers for their kind invitation.}\\
Universit\'e de Paris-Sud
\\
pascal.auscher@math.u-psud.fr
\and
Thierry Coulhon \\
Universit\'e de Cergy-Pontoise\\
thierry.coulhon@math.u-cergy.fr }

\date{April 11, 2005}

\maketitle

\abstract{We study the validity of the $L^p$ inequality for the  Riesz transform    when $p>2$ and of its reverse inequality when $p<2$  on  complete Riemannian manifolds  under the doubling property and some Poincar\'e inequalities.\\

{\bf MSC numbers 2000:} 58J35, 42B20
}

\tableofcontents

 \section*{Introduction}

Let $M$ be a non-compact complete  Riemannian manifold.  Denote by  $\mu$ the  Riemannian measure, and by $\nabla$  the
Riemannian gradient. 
Denote by $|.|$ the length in the tangent space,
and by $\|.\|_p$ the norm in $L^p(M,\mu)$, $1\le p\le \infty$. One defines $\Delta$, the
Laplace-Beltrami operator, as a self-adjoint positive operator on $L^2(M,\mu)$ by the formal
integration by parts 
\begin{equation*}
(\Delta f,f)= \||\nabla f|\|_2^2
\end{equation*} 
for all 
$f\in\ccc^\infty_0(M)$, and its positive self-adjoint square root $\Delta^{1/2}$ by
\begin{equation*}
(\Delta f,f)=\|\Delta^{1/2}f\|_2^2.
\end{equation*}  As a consequence,
$$
\norm{|\nabla f|}_2^2=\|\Delta^{1/2}f\|_2^2.  \eqno(E_{2})
$$
To identify the spaces defined by (completion with respect to) the seminorms
$\norm{|\nabla f|}_p$ and $\norm{\Delta^{1/2} f}_p$  on $\ccc_0^\infty(M)$ for some $p\in
(1,\infty)$, it is enough to prove that there exist $0<c_p\le
C_p<\infty$ such that for all $f\in
\ccc_0^\infty(M)$
$$ 
c_p\norm{\Delta^{1/2} f}_p\le \norm{|\nabla f|}_p\leq C_p\norm{\Delta^{1/2} f}_p. \eqno(E_{p})
$$
This equivalence splits into two inequalities of different nature. The right-hand inequality may be reformulated by saying that 
 the Riesz transform $\nabla \Delta^{-1/2}$
   is bounded from $L^p(M,\mu)$ to the   space of  $L^p$ vector fields\footnote{In the case where $M$ has finite measure,
 one has to replace $L^p(M)$ by the subspace
$L^p_0(M)$ of  func\-tions  with mean zero; this modification will be implicit in what follows.},   in other
words
$$\norm{|\nabla \Delta^{-1/2}f|}_p \le C_p
 \|f\|_p\, . 
\eqno(R_p)$$
 The left hand inequality is what we call the reverse inequality
$$
\norm{\Delta^{1/2} f}_p\leq C_p\norm{|\nabla f|}_p.
\eqno(RR_p)$$
 
It is well-known (see \cite{Ba0}, Section 4,
or
\cite{CDfull}, Section 2.1) that   $(R_p)$  
 implies  $(RR_{p'})$
where $p'$ is the conjugate exponent of $p$ but the converse is not clear (in fact, it is false, see below). We mention a partial converse which we shall use and prove in the sequel.

\begin{lem}\label{reciproque} The conjunction of $(RR_{p'})$ and $(\Pi_{p})$ implies $(R_{p})$.
\end{lem}

Here, $(\Pi_{p})$ is the inequality describing the boundedness  on $L^pT^*M$ of the orthogonal projector $d\Delta^{-1}\delta$ of $1$-forms onto exact forms. Namely, for all $\omega \in \ccc_{0}^\infty(T^*M)$,
$$
\norm{|d\Delta^{-1}\delta \omega|}_{p}\le C_{p}\norm \omega _{p}, \eqno(\Pi_{p})
$$
where $d$ is the exterior derivative and $\delta$ its formal adjoint. 
\bigskip

The question is to find which geometrical properties on $M$ insure each of these inequalities, and in the end $(E_{p})$ for a range of $p$'s.

We first recall the result of \cite{CD} which deals with $(R_p)$ for $1<p<2$.
Denote by  $B(x,r)$ the open ball of radius $r>0$ and center $x\in M$, and by $V(x,r)$ its
measure $\mu(B(x,r))$. One says  that $M$ satisfies the doubling
property if there exists $C>0$ such that, for all $x\in M$ and $r>0$,
$$
V(x,2r)\le C\,V(x,r).\eqno (D)
$$
Let  $p_t(x,y)$, $t>0$, $x,y\in M$ be the heat kernel of $M$, that is the kernel of the heat semigroup $e^{-t\Delta}$.

\begin{theo}[\cite{CD}]\label{CD} Let $M$ be a complete non-compact Riemannian manifold satisfying $(D)$. Assume that for all $x\in M$, $t>0$ and some constant $C>0$,
$$p_t(x,x)\le \frac{C}{ V(x,\sqrt{t})}. \eqno(DU\!E)$$
Then  $(R_{p})$ holds for $1<p<2$, hence $(RR_{p})$ for $2<p<\infty$.
\end{theo}

It is also shown in
\cite{CD} that  the Riesz transform is unbounded on $L^p$  for every $p>2$ on the manifold consisting of two copies of the
Euclidean plane glued smoothly along their unit circles, although it satisfies
$(D)$ and $(DU\!E)$.  

A stronger
assumption is therefore required  to obtain $(R_{p})$ when $p>2$. 

It is natural to assume in addition the Poincar\'e inequalities,  although it is known that they are not sufficient for  $(R_{p})$ to hold for all $p>2$ (\cite{LHQ}, \cite{CouLi}),  nor necessary for
$(R_{p})$ to hold for some $p>2$ (\cite{CCH}).
One says that $M$ satisfies the (scaled) Poincar\'e inequalities    $(P_{2})$ if there exists $C>0$ such that, for every ball $B=B(x,r)$, $x\in M$, $r>0$, and
every
$f$ with
$f,
\nabla f$ locally in $L^2$,
 $$
\int_{B}|f-f_{B}|^2\,d\mu\le Cr^2\int_{B}
|\nabla f|^2\,d\mu, \eqno (P_{2})
$$
where $f_E$ denotes the mean of $f$ on the set $E$. 

Even under $(D)$ and $(P_2)$ alone,  it is not clear that $(R_{p})$ holds for some $p>2$ because of the following result proved in \cite{ACDH} which tells us that the semigroup should 
have some boundedness properties (it is also shown there that this is the same as some 
$L^p$ estimates of the gradient of the heat kernel).

 \begin{theo}\label{nsc} Let $M$ be a complete non-compact Riemannian manifold satisfying $(D)$
and
$(P_2)$. Let $p_0 \in (2,\infty]$. The following assertions are
equivalent:
\begin{enumerate}
\item
  For all  $p \in (2,p_0)$, there exists $C_p$ such that  for  all $t>0$
 $$\||\nabla e^{- t\Delta}|\|_{p\to p} \le  \frac{C_p}{\sqrt t}.
$$
\item $(R_{p})$ holds for $p\in (2,p_0)$. 
\end{enumerate}
\end{theo}

 Our main result states that, in the situation of Theorem \ref{nsc}, there always exists a $p_0=2+\ep>2$ such that condition 2 is satisfied. 
 
  \begin{theo}\label{Rp>2} Let $M$ be a complete non-compact Riemannian manifold satisfying $(D)$ and $(P_{2})$.  Then there is   $\ep>0$ such that    $(R_p)$ holds for $2<p<2+\ep$. \end{theo}
  
 Our proof does not rely on Theorem \ref{nsc} and in fact we shall add a list of  assertions equivalent to condition 2, one of them being easier to check. But in view of Theorem \ref{nsc}, this also says that there is an automatic improvement of $L^p$ estimates for the gradient of the semigroup, which is reminiscent (and, as we shall see, equivalent) to the self-improvement  ``\`a la Meyers'' of Sobolev $W^{1,p}$ estimates for weak solutions of elliptic equations \cite{Ms}.

It is well-known (see \cite{S2}, \cite{parma}) that the conjonction of $(D)$ and $(P_{2})$ is equivalent to the full Li-Yau type estimate
$$
 \frac{c}{ V(y,\sqrt{t})}\exp\left(-C  \frac{d^2(x,y)}{
t}\right)\le p_t(x,y)\le \frac{C}{ V(
y,\sqrt{t})}\exp\left(-c  \frac{d^2(x,y)}{
t}\right), \eqno (LY)
$$
for all $x,y\in M$, $t>0$ and some constants $C,c>0$. Hence, $(D)$ and $(P_{2})$ imply $(D)$ and $(DU\!E)$. Therefore combining Theorems \ref{CD} and \ref{Rp>2},  we obtain

\begin{cor}\label{main} Let $M$ be a complete non-compact Riemannian manifold satisfying $(D)$ and $(P_{2})$. Then there is   $p_0 \in (2,\infty)$ such that  $(E_p)$ holds  when  $p'_{0}<p<p_0$.
\end{cor}

A crucial step towards Theorem  \ref{Rp>2} consists in giving a sufficient condition for the reverse inequality $(RR_p)$ for $1<p<2$ in terms of
the $L^p$ version of $(P_2)$. Let $1\le p<\infty$. One says that $M$ satisfies   $(P_{p})$ if there exists $C>0$ such that, for every ball $B=B(x,r)$ and
every
$f$ with
$f,
\nabla f$ locally $p$-integrable,
 $$
\int_{B}|f-f_{B}|^p\,d\mu\le Cr^p\int_{B}
|\nabla f|^p\,d\mu. \eqno (P_{p})
$$ 
It is known that $(P_{p})$ implies $(P_{q})$ when $p<q$ (see for instance \cite{HaK}).  Thus the set of $p$'s such that $(P_p)$ holds  is,  if it is not empty, an interval unbounded on the right.  A recent deep result asserts in a general context of metric measured spaces that this interval is open in $[1,+\infty[$. In our case, it states as follows. 

\begin{lem}[\cite{Ke}]\label{selfimprovement} Let $M$ be a complete non-compact Riemannian manifold satisfying $(D)$. Assume $p>1$. Then $(P_{p})$ self-improves to $(P_{p-\ep})$ for some $\ep>0$.
\end{lem}

We shall prove

\begin{theo}\label{poincareq} Let $M$ be a complete non-compact Riemannian manifold satisfying $(D)$ and $(P_{q})$  for some $q \in [1,2]$.  Then    $(RR_p)$ holds for $q<p<2$.  If $q=1$, there is a weak type $(1,1)$ estimate.
\end{theo}

Define $q_0=\inf\{p\in [1,2]; (P_p) \mbox{ holds}\}$. Note that if $(P_p)$ holds for some $p\in (1,2]$, then $q_0<p$ according to 
Lemma \ref{selfimprovement}.
As a consequence of Theorem \ref{poincareq} and Lemma \ref{selfimprovement}, if $q_0<2$, 
that is to say if $(P_2)$ holds,
$(RR_p)$ holds for $p\in (q_0,2]$.

As a corollary of Theorems  \ref{CD}, \ref{Rp>2} and \ref{poincareq} we obtain for instance

\begin{cor}\label{mainc} Let $M$ be a complete non-compact Riemannian manifold satisfying $(D)$ and $(P_{1})$. Then   $(E_p)$ holds  when  $1<p<2+\ep$ for some $\ep>0$. 
\end{cor}

One may observe that our proofs do not use completeness in itself, but rather stochastic completeness,
that is the property
\begin{equation}
\int_Mp_t(x,y)\,d\mu(y)=1,\label{stoco}
\end{equation}
for all $x\in M$ and $t>0$,
which does hold for complete manifolds satisfying $(D)$ (see \cite{Grigostoc}),
but also for instance  for  conical  manifolds  with closed basis (see \cite{LHQ}).

Note that the class of manifolds satisfying $(D)$ and $(P_1)$ (therefore also $(P_2)$)
contains all complete manifolds that are quasi-isometric to a manifold with non-negative Ricci curvature
(see \cite{parma}).

It is proved in \cite{CDfull} that for any $q\in (1,2)$, there exists a complete Riemannian manifold with $(D)$ such that $(RR_{p})$ fails for all $1<p<q$ \footnote{\,We remark that the positive result in \cite{CDfull} concerning $(RR_{p})$, namely Theorem 6.1, has a gap, since it depends from another result in the same paper, Proposition 5.4, which has a mistake in the argument. The mistake is located in the last line of p. 1744 where it is said that the (usual) Calder\'on-Zygmund decomposition preserves exact forms. This is exactly the obstacle that we get around in Section 1 with a modified Calder\'on-Zygmund decomposition and it is not clear that the same ideas can be employed  under the assumption taken in \cite{CDfull}.}.
The point is that  there are manifolds satisfying a $L^2$ Sobolev 
inequality at infinity associated with
a certain dimension, but, for $p$ close to $1$, only a $L^p$ Sobolev inequality  
associated with a much  lower dimension,
and, for $p=1$, a trivial isoperimetric inequality,
whereas $(RR_p)$ would impose a tighter connection between $L^2$ and $L^p$ 
Sobolev inequalities. In other words, $(RR_p)$ imposes that
the heat kernel dimension and the isoperimetric dimension cannot differ too 
much.

It has been proved by Li Hong-Quan in \cite{LHQ}  that, on conical manifolds with closed basis,
$(R_p)$ holds if and only if $1<p<p_0$, where the threshold $p_0>2$ depends on the $\lambda_1$ of the basis. Now, all these manifolds satisfy $(P_2)$ (see \cite{CouLi}) and one can see that they even satisfy  $(P_1)$ by using the methods in \cite{GS}.   In particular, there is no hope that the assumptions
of Corollary  \ref{mainc} suffice for $(R_p)$ to hold for all $p>2$.

In view of Corollary \ref{mainc}, this also shows that, as we mentioned above,  $(RR_p)$  does not imply 
  $(R_{p'})$, even in the class of manifolds with doubling, in the range $1<p<2$.

 Let us summarize the situation for (stochastically) complete Riemannian manifolds, satisfying $(D)$, going from weakest to strongest hypotheses.
\begin{enumerate}
  \item It is known that $(R_{p})$ may be false for $2<p$ and that $(RR_{p})$ may be false for 
  $1<p<2$. What can be said about the other cases, that is $(R_{p})$ for $p<2$ and $(RR_{p})$ for $p>2$?
  \item Assume  $(DUE)$. Then $(R_{p})$ holds for $1<p \le 2$, $(RR_{p})$ for $p\ge 2$ and $(R_{p})$ may be false for all $p>2$. What can be said about $(RR_{p})$ for $p<2$?
  \item  Assume $(P_{2})$. Then $(R_{p})$ holds for $1<p<p_{0}$ with some $p_{0}>2$, $(RR_{p})$ for $q_{0}<p<\infty$ with some $q_{0}<2$. 
  Can one give estimates on $p_{0}$ and $q_{0}$?  
   \item  Assume $(P_{1})$. Then $(R_{p})$ holds for $1<p<p_{0}$ with some $p_{0}>2$, $(RR_{p})$ for $1<p<\infty$.
  Can one give estimates on $p_{0}$?
  \end{enumerate}


The proof of Theorem \ref{poincareq} in Section \ref{secRRp} uses methods of the first author in \cite{Au0} adapted to the present situation and in particular a Calder\'on-Zygmund lemma for Sobolev functions, which allows us to do a Marcinkiewicz type interpolation.

As said, we do not rely  on Theorem \ref{nsc} to prove Theorem \ref{Rp>2}. Instead, we use ideas of Shen in \cite{Shen} developed for elliptic operators on Euclidean space and extend them to the class of manifolds we consider.  This yields a new characterization of the $L^p$ boundedness of Riesz transforms for $p>2$ (with a restriction that $p$ should be close to 2) in terms of local and scale invariant estimates on harmonic functions (Theorem \ref{newCNS}) which are more tractable in practice.  
 In passing, we show that this is also equivalent to the $L^p$ boundedness of $d\Delta^{-1}\delta $.  Actually the main tool in \cite{Shen}  is a theorem (Theorem 3.1)  for boundedness of operators with no kernels which is essentially similar to Theorem 2.1 in \cite{ACDH}.

\section{Reverse inequalities $(RR_{p})$ for $p<2$}\label{secRRp}

In this section, we prove  Theorem \ref{poincareq}. We assume that  $(D)$ and  $(P_{q})$ hold for a $1\le q<2$ and prove 
$(RR_{p})$ for $q<p<2$.

We first establish a Calder\'on-Zygmund lemma for Sobolev functions. Next, we apply this lemma to establish  the  preliminary  weak type estimate
\begin{equation}
\label{weaktypeq}
\norm{\Delta^{1/2} f}_{q,\infty} \leq C_q\norm{|\nabla f|}_q, \ \forall\,
f\in\ccc^\infty_0(M).
\end{equation}
Finally,  we proceed via an interpolation argument.

\subsection{A Calder\'on-Zygmund lemma for Sobolev functions}

We present here in the Riemannian context a result first proved by one of us \cite{Au0} in the Euclidean setting with Lebesgue measure  (see also  the extension to weighted  Lebesgue measure in \cite{AusMar}).

\begin{pro}\label{lemmaCZD-manifold} Let $M$ be a complete non-compact Riemannian manifold satisfying $(D)$. Let $1\le q< \infty$ and assume that $(P_{q})$ holds. Let  $f\in \ccc^\infty_0(M)$\,\footnote{Of course, $f$ can be taken more general than this.} be such that
$\norm{|\nabla f|}_{q} <\infty$. Let $\alpha>0$.\footnote{If $\mu(M)<\infty$, one has to assume $\mu(M)\alpha^q >C\int |\nabla f|^q\, d\mu$ for some constant $C$ depending only on $(D)$.} Then, one can find a
collection of balls $B_i$,  $C^1$ functions  $b_i$ and a (almost everywhere) Lipschitz function $g$  such that the following properties hold:
\begin{equation}\label{eqcsds1}
f= g+\sum_i b_i,  \end{equation}
\begin{equation}\label{eqcsds2}
|\nabla g(x)| \le C\alpha, \quad \text{for}\  \mu-\text{a.e.} \ x\in M,
\end{equation}
\begin{equation}\label{eqcsds3}
\supp b_i \subset B_{i}\ \text{and} \ \int_{B_i} |\nabla b_i|^q\, d\mu \le
C\alpha^q \mu(B_i), \end{equation}
 \begin{equation}\label{eqcsds4}
\sum_i \mu(B_i) \le C\alpha^{-q} \int |\nabla f|^q\, d\mu ,
\end{equation}
\begin{equation}\label{eqcsds5}
\sum_i {\bf 1}_{B_i} \le N, \end{equation}
where $C$ and
 $N$ only depend on $q$ and on the  constant in $(D)$. 
\end{pro}

\paragraph{Proof:} Let  $f\in \ccc^\infty_0(M)$ and $\alpha>0$.
Consider $\Omega= \{x \in M; {\cal M}(|\nabla f|^q)(x) >\alpha^q\}$, where
${\cal M}$ is the uncentered maximal operator over balls of $M$. If $\Omega$ is empty, then set $g=f$, $b_i=0$; (\ref{eqcsds2}) is satisfied thanks to  Lebesgue differentiation theorem. Otherwise, the maximal theorem gives us
\begin{equation}\label{max}
\mu(\Omega) \le C\alpha^{-q} \int  |\nabla f|^q\, d\mu.  \end{equation}
Let $F$ be the complement of $\Omega$. Again by the Lebesgue differentiation
theorem, $|\nabla f| \le \alpha$  $\mu$-almost everywhere on $F$.  Since $\Omega$ is open, let $(\underline B_i)$ be a Whitney decomposition of $\Omega$. That is,
$\Omega$ is the union of the $\underline B_i$'s, and there are constants $C_{2}>C_{1}>1$ depending only on the metric such that 
the balls $B_i= C_{1}\underline B_{i}$ are contained in  $\Omega$ and have the bounded overlap
property, but each ball $\overline B_i= C_{2}\underline B_{i}$ intersects $F$ (see \cite{CW}).
As usual, $C B$ is the ball co-centered with $B$ with radius $C r(B)$.
Condition 
\eqref{eqcsds5} is nothing but the bounded overlap property and  \eqref{eqcsds4} follows from  \eqref{eqcsds5} and \eqref{max}. Furthermore,  $\overline{B_i} \cap F\ne \emptyset$ and the doubling property imply
$$
\int_{B_i} |\nabla f|^q \, d\mu \le \int_{\overline{B_i}} |\nabla f|^q \, d\mu \le \alpha^q \mu(\overline{B_i}) \le C\alpha^q \mu(B_i).
$$
Let us now define the
functions $b_i$. Let $(\calX_i)$ be a partition of unity of $\Omega$
subordinated to the covering $(\underline B_i)$ so that for each $i$, $\calX_i$ is a $C^1$
function supported in $B_i$ with $
 \|\nabla \calX_i\|_\infty \le \frac{C}{r_i}$, $r_i=r(B_{i})$. Set
$$
b_i = (f-f_{B_{i}})\calX_i.
$$
It is clear that $b_i$ is supported in $B_i$. Let us estimate $\int_{B_i}
|\nabla b_i|^q\, d\mu$.
Since $$\nabla \left(\,(f-f_{B_{i}})\calX_i\right) = \calX_i \nabla f +
(f-f_{B_{i}})\nabla\calX_i,$$
we have  by the $L^q$ Poincar\'e inequality and the above estimate on $\nabla f$ that $$\int_{B_i} |\nabla \left(\,(f-f_{B_{i}})\calX_i\right)|^q\, d\mu  \le
C\alpha^q \mu(B_i).$$ Thus \eqref{eqcsds3} is proved.

Set $g=f-\sum b_{i}$. Then $g$ is defined $\mu$-almost everywhere since the sum is locally finite on $\Omega$ and vanishes on $F$, and $g$ is also defined in the sense of distributions on $M$ (not just $\Omega$ which is trivial:  in fact the argument shows that  $g$ is a locally integrable function on $M$). For the latter claim, if $\varphi\in \ccc^\infty_0(M)$, we observe that for $x$ in the support of $b_i$, we have $d(x,F)\ge r_i$, so that  
$$
\int \sum_{i } |b_{i}|| \varphi|\, d\mu \le \bigg(\int \sum_{i}\frac{|b_{i}|}{r_{i}}\, d\mu\bigg) \ \sup_{x\in M} (d(x,F) |\varphi(x)|). 
$$
 By H\"older inequality and the Poincar\'e $L^q$ inequality, 
$$
\int \frac{|b_{i}|}{r_{i}}\, d\mu \le (\mu(B_i))^{1/q'} \left (\int_{B_i} |\nabla f|^q \, d\mu\right)^{1/q} \le C\alpha \mu(B_i).$$
Hence
$$
\int \sum_{i } |b_{i}|| \varphi|\, d\mu \le C\alpha\, \mu(\Omega) \ \sup_{x\in M} (d(x,F) |\varphi(x)|), 
$$
which proves the claim.

It remains to prove \eqref{eqcsds2}.  Note that
 $\sum_i \calX_i(x)= 1$ and  $\sum_i  \nabla\calX_i(x)=0$ for $x
\in \Omega$.  It follows that
\begin{eqnarray*} \nabla g &= &\nabla f - \sum_{i} \nabla b_{i}\\
&=&\nabla f -(\sum_i \calX_i) \nabla f -
\sum_i(f-f_{B_{i}})\nabla\calX_i\\
&=& (\nabla f) {\bf 1}_F+\sum_i f_{B_{i}} \ \nabla\calX_i.
\end{eqnarray*}

Note that by the definition of $F$,
$|(\nabla f) {\bf 1}_F| \le \alpha$.
We claim that a similar estimate holds for
 $h= \sum_i f_{B_{i}} \ \nabla\calX_i$, that is $|h(x)| \le C\alpha$ for all $x\in M$ for some constant $C$ independent of $x$. Note that this
sum vanishes on $F$ and is locally finite on $\Omega$.  
Fix now $x \in \Omega$. Let $B_j$ be some Whitney ball
containing $x$ and let $I_x$ be the set of indices $i$ such that
$x \in B_i$. We know that $\sharp I_x \le N$.  Also for $i \in
I_x$ we have that $C^{-1}r_i \le  r_j \le Cr_i$ where the
constant $C$ depends only on doubling (see \cite[Chapter I,  3]{St1} for the Euclidean case).  We also
have  $|f_{B_{i}}-f_{B_{j}}|\le Cr_j\alpha$. 
Indeed, one has $B_i \subset AB_j$ with $A=2C+1$, so that  by the
Poincar\'e $L^q$ inequality one obtains
\begin{align*}|f_{B_{i}}- f_{AB_{j}}| &\le \frac{1}{\mu(B_i)}\int_{B_i}|f- f_{AB_{j}}| \\
&\le  \frac{C}{\mu(B_{j})}\int_{AB_{j}}|f- f_{AB_{j}}|
\\
&\le C  Ar_j ((|\nabla f|^q)_{AB_{j}})^{1/q}
\\
&\le
CAr_j\alpha
\end{align*}
and similarly for $|f_{AB_j} - f_{B_j}|$.
Hence,
$$
|h(x)| = \left|\sum_{i \in I_x} (f_{B_{i}}-f_{B_{j}})
\nabla\calX_i(x)\right| \le C \sum_{i \in I_x} |f_{B_{i}}-f_{B_{j}}
|r_i^{-1} \le C N \alpha.$$
This proves \eqref{eqcsds2}, and finishes the proof of  Proposition \ref{lemmaCZD-manifold}.

\paragraph{Remarks} 1) It follows from the construction that $\sum \nabla b_i \in L^q$ with 
norm bounded by $C\norm{|\nabla f|}_q$, hence $\norm{|\nabla g|}_q \le (C+1) \norm{|\nabla f|}_q$.

2) $g$ is equal almost everywhere to a  Lipschitz function on $M$ and $|g(x)-g(y)| \le C\alpha d(x,y)$ almost everywhere. The point is that the Lipschitz constant is controlled by $\alpha$. This can be shown by similar arguments as for obtaining \eqref{eqcsds1}. Alternatively, once \eqref{eqcsds1} is proved, one can show that $g, \nabla g$ satisfy the $q$-Poincar\'e inequality on arbitrary balls by using the definition of $g$ as $f-\sum b_{i}$ since $f$ and each $b_{i}$ do. At this point, we invoke Theorem 3.2 in \cite{HaK} and the $L^\infty$ bound on $|\nabla g|$ to conclude. 

3)  Observe that 
$
g=f\textbf{1}_F + \sum f_{B_i} \calX_i
$
so that is contains in particular the fact that $f$ is equal almost everywhere to a Lipschitz function on $F$.
Hence, $g$ is some sort of Whitney extension of the restriction of $f$ to $F$ where averages of $f$ on $B_i$ (since $f$ was already defined on the complement of $F$) replace evaluation at some point inside $F$ at distance
$Cr_i$ to $B_i$.

\subsection{A weak type estimate}

Assume $(P_q)$ for some $q\in[1,2)$. 
Let $f\in 
\ccc_0^\infty(M)$. We wish to establish the    estimate
\begin{equation}\label{eq21}
\mu\left(\left\{x \in M; |\Delta^{1/2}f(x)| >\alpha\right\}\right) \le 
\frac{C}{\alpha^q}\int |\nabla f|^q\, d\mu,
\end{equation}
for all $\alpha>0$.  We use the following resolution of $\Delta^{1/2}$:
$$
\Delta^{1/2}f= c\int_0^\infty \Delta e^{-t \Delta}  f \, \frac{dt}{\sqrt t}
$$
where $c=\pi^{-1/2}$ is forgotten from now on. It suffices to obtain the 
result for the truncated integrals $\int_\ep^R\ldots$ with bounds independent of $\ep,R$, 
and then to let
$\ep\downarrow 0$ and $R\uparrow \infty$. 
For the truncated integrals, all the calculations are justified. We henceforth assume that $\Delta^{1/2}$ is replaced by  one of the truncations above but we keep writing $\Delta^{1/2}$ and the limits of the integral as $0, \infty$ to keep the notation   simple.
 
 Assume first  $\mu(M)=\infty$. Apply the Cal\-de\-r\'on-Zygmund decomposition of Proposition \ref{lemmaCZD-manifold} 
to $f$ at height $\alpha$ with exponant $q$ and write $f=g+\sum_i b_i$.

Since $g$ and $b_i$ are no longer $C^\infty_0(M)$ we have to give a meaning to 
$\Delta^{1/2} g$ and $\Delta^{1/2} b_i$. Since $\Delta^{1/2}$ is replaced by approximations, it suffices to 
define $\Delta e^{ -t\Delta}g$ and $\Delta e^{-t\Delta} b_i$ for $t>0$. Since $(D)$ and $(P_q)$ imply $(D)$
and $(P_2)$, we have the Gaussian upper bounds for the kernel of $e^{-t\Delta}$ and by analyticity for the kernel of $t\Delta e^{-t\Delta}$. Since $b_i$ has support in a ball and is integrable (see the proof of Proposition \ref{lemmaCZD-manifold}) $\Delta e^{-t\Delta} b_i(x)$ is defined by the convergent integral
$\int_M \partial_t p_t(x,y) b_i(y)\, d\mu(y)$.

As for $g$, we know it  equals almost everywhere a Lipschitz function with Lipschitz constant bounded by $C\alpha$  (see Remarks 1 and 2 at the end of Section 1.1).
We fix any point $z$ where $g(z)$ exists and we have that 
 $\int_M  \partial_t p_t(x,y) g(y)\, d\mu(y)$ is a smooth function bounded by $C\alpha t^{-1} (d(x,z)+ t^{1/2})$ (we use the fact
 that $\int_M \partial_tp_t(x,y)\, d\mu(y)=0$). We take this as our definition of $\Delta e^{-t\Delta} g(x)$. 
 
 Next, we prove
 \begin{equation*}\label{eq22}
\mu\left\{x \in M; |\Delta^{1/2}g(x)| >\frac \alpha 3 \right\}  \le 
\frac{C}{\alpha^q}\int_M |\nabla f|^q \, d\mu.
\end{equation*}
Since
\begin{equation*}
\mu\left\{x \in M; |\Delta^{1/2}g(x)| >\frac \alpha 3 \right\} \le 
\frac{9}{\alpha^2}\int_M
|\Delta^{1/2} g|^2\, d\mu,
\end{equation*}
 it remains to justify
 \begin{equation} \label{eq:L1/2g}
\int_M 
|\Delta^{1/2} g|^2\, d\mu \le \int_M |\nabla g|^2\, d\mu.
\end{equation}
Indeed, once this is done, we conclude by using  $\int_M |\nabla g|^2\, d\mu \le C\alpha^{2-q}
 \int_M |\nabla f|^q\, d\mu$ which follows from $\norm{|\nabla g|}_q \le C\norm{|\nabla f|}_q$ and \eqref{eqcsds2} since $q<2$.

Note that \eqref{eq:L1/2g} (since we have replaced $\Delta^{1/2}$ by truncations) would be valid  if 
$g$ were in $C_0^\infty(M)$. 
  For $\varphi\in C_0^\infty(M)$, we have by Fubini's theorem
 \begin{align*}\int_M \Delta e^{-t\Delta} g(x) \varphi(x) \, d\mu(x)
 &= 
 \int_M g(y)   \Delta e^{-t\Delta}	\varphi (y) \, d\mu(y) 
 \\
 &
 =
 \lim_{r\to +\infty}  \int_M \eta_r(y) g(y)   \Delta e^{-t\Delta} \varphi(y) \, d\mu(y).
 \end{align*}
 Here $\eta_r$ is a smooth cut-off which is bounded by 1 on $M$, equal to 1 on a ball   $B_r$  of radius $r$, 0 outside the ball $2B_r$, and with $\norm{|\nabla \eta_r|}_\infty \le C/r$.
 By Stokes theorem, the last integral is equal to 
 $$
 \int_M \eta_r \nabla g \cdot \nabla e^{-t\Delta} \varphi \, d\mu + \int_M g\nabla \eta_r \cdot \nabla e^{-t\Delta} \varphi \, d\mu.$$
 Under our assumptions,  we have the weighted $L^2$ estimate from \cite{G1}  (see also \cite{CD}):
for some $\gamma>0$ and all $y\in M, t>0$, 
\begin{equation}\label{gradx}
\int_M | \nabla_x\, p_t(x,y)|^2  e^{\gamma\frac{ d^2(x,y)}{t}} \, d\mu(x) \le \frac C{ t
\,V(y,\sqrt{t})}
\end{equation}
where $\nabla_x$ means that the gradient is taken with respect to the $x$ variable. 
Given the fact that $\nabla g$ is square integrable  and $g$ is Lipschitz, it is not difficult to pass to the limit as $r\to \infty$ and to conclude that 
$$
\int_M \Delta e^{-t\Delta} g\,  \varphi \, d\mu=\int_M \nabla g \cdot \nabla e^{-t\Delta} \varphi \, d\mu.
$$
Thus, we obtain (again, $\Delta^{1/2}$ is replaced by truncated integrals)
$$
\langle \Delta^{1/2} g, \varphi\rangle = 
\langle \nabla g , \nabla \Delta^{-1/2} \varphi\rangle,
$$
so that  a duality argument from the equality $(E_2)$ (or rather its approximation)  yields 
\eqref{eq:L1/2g}.

To compute $\Delta^{1/2}b_i$, let $r_i=2^k$ if $2^k \le  r(B_i) < 2^{k+1}$ and
set 
$T_i= \int_0^{r_i^2} \Delta e^{-t \Delta}\, \frac{dt}{\sqrt t}$ and
 $U_i= \int_{r_i^2}^\infty \Delta e^{-t \Delta} \, \frac{dt}{\sqrt t}$. It is enough to 
estimate $A=\mu\{x \in M; |\sum_i T_ib_i(x)| >\alpha/3\}$ and 
$B=\mu\{x \in M; |\sum_i U_ib_i(x)| >\alpha/3\}$. 

First
$$A \le \mu(\cup_i 4B_i) + \mu\left(\left\{x \in M \setminus \cup_i  4B_i ; 
\left|\sum_i T_ib_i(x)\right| >\frac \alpha 3\right\}\right),$$
and by \eqref{eqcsds4} and $(D)$, $\mu(\cup_i 4B_i) \le \frac{C}{\alpha^q}\int_M  |\nabla f|^q \, d\mu$.

For the other term, we have
$$
\mu\left(\left\{x \in M \setminus \cup_i 4B_i ; \left|\sum_i
T_ib_i(x)\right| >\frac \alpha 3\right\}\right)
\le
\frac{C}{\alpha^2}\int_M \left|\sum_i h_i\right|^2\, d\mu$$ with $h_i = {\bf
1}_{(4B_i)^c}|T_ib_i|$. To estimate the $L^2$ norm, we follow ideas in \cite{BK1, HM} and dualize against
$u\in L^{2}(M,\mu)$ with $\|u\|_{2}=1$ and write
$$
\int_M  |u| \sum_i h_i \, d\mu = \sum_i\sum_{j=2}^\infty A_{ij}
$$  
where 
$$
A_{ij}= \int_{C_j(B_i)} |T_ib_i||u|\, d\mu
$$   
with $C_j(B_i) = 2^{j+1}B_{i} \setminus 2^jB_{i}$. By Minkowski integral inequality \begin{equation*}\| T_ib_i \|_{L^2(C_j(B_i))} \le   \int_0^{r_i^2} \|\Delta e^{-t \Delta}  b_i \|_{L^2(C_j(B_i))} \, \frac{dt} {\sqrt t} 
\end{equation*}
and by the  Gaussian upper bounds for the kernel of 
 $\Delta e^{-t \Delta} $ (see above),
$$
|\Delta e^{-t \Delta} b_i (x)| \le \int_{M} \frac {C}{tV(y,\sqrt t)} e^{-\frac{cd^2(x,y)}{t}} \, |b_{i}(y)|\, d\mu(y).
$$
Now, $y$ is in the support of $b_{i}$, that is $B_{i}$, and $x \in C_{j}(B_{i})$, hence, one may replace $d(x,y)$ by $2^jr_i$ in the Gaussian term since $r_i \sim r(B_{i})$. Also, if $y_{i}$ denotes the center of $B_{i}$, write
$$
\frac{V(y_{i},\sqrt t)}{V(y, \sqrt t)} = \frac{V(y_{i},\sqrt t)}{V(y_{i}, r_{i})}  \frac{V(y_{i},r_{i})}{V(y, r_{i})} \frac{V(y,r_{i})}{V(y, \sqrt t)}.
$$
By $(D)$ and $\frac {V(z,r)}{V(z,s)} \le c ( \frac r s )^\beta$ for $r>s$, as $t \le r_{i}^2$, we have
$$
\frac{V(y_{i},\sqrt t)}{V(y, \sqrt t)} \le c \left ( \frac{r_{i}}{\sqrt t }\right)^\beta.
$$
Using this estimate, 
$
\int_{B_i} |b_{i}|\, d\mu \le C \mu(B_{i}) r_{i} \alpha 
$
and $ \mu(B_{i}) \sim V(y_{i}, r_{i})$, we obtain
\begin{align*}|\Delta e^{-t \Delta} b_i (x)| &\le \frac {C}{tV(y_{i},\sqrt t)}   \left ( \frac{r_{i}}{\sqrt t }\right)^\beta e^{-\frac{c 4^j r_{i}^2}{t}} \int_{B_i} |b_{i}|\, d\mu
\\
&\le  \frac {C r_{i}}{t}   \left ( \frac{r_{i}}{\sqrt t }\right)^{2\beta} e^{-\frac{c 4^j r_{i}^2}{t}} \alpha.
\end{align*}
Thus, 
$$
\|\Delta e^{-t \Delta}  b_i \|_{L^2(C_j(B_i))} \le   \frac {C r_{i}}{t}   \left ( \frac{r_{i}}{\sqrt t }\right)^{2\beta} e^{-\frac{c 4^j r_{i}^2}{t}}  (\mu(2^{j+1}B_{i}))^{1/2} \alpha.
$$
Plugging this estimate  inside the integral, we obtain
$$
 \| T_ib_i \|_{L^2(C_j(B_i))} \le   Ce^{-c4^j}(\mu(2^{j+1}B_{i})^{1/2}\alpha
$$
for some $C,c>0$.  




Now remark that for any $y \in B_i$ and any $j\ge 2$, 
$$
\left( \int_{C_j(B_i)} |u|^{2}\, d\mu\right)^{1/2} \le \left( \int_{2^{j+1}B_i} |u|^{2}\, d\mu \right)^{1/2} \le  \mu(2^{j+1}B_i)^{1/2}
\big({\cal M}(|u|^{2})(y)\big)^{1/2}.
$$
Applying H\"older inequality and doubling,    one obtains
$$
A_{ij} \le  C\alpha 2^{j\beta} e^{-c4^j} \mu(B_{i}) \big({\cal M}(|u|^{2})(y)\big)^{1/2}.
$$
Averaging over $y\in B_i$ yields
$$
A_{ij} \le C\alpha 2^{j\beta} e^{-c4^j} \int_{B_i} \big({\cal M}(|u|^{2})\big)^{1/2}\, d\mu.$$
Summing over $j\ge 2$ and $i$, we have
$$
\int_M | u| \sum_i h_i\, d\mu \le C \alpha \int_M  \sum_i {\bf 1}_{B_i} \big({\cal M}(|u|^{2})\big)^{1/2}\, d\mu.$$
Using finite overlap \eqref{eqcsds5} of the balls $B_i$ and Kolmogorov's inequality, one obtains
$$
\int_M | u| \sum_i h_i\, d\mu  \le C'N\alpha\mu \big( \cup_i B_i \big)^{1/2} \||u|^2\|_1^{1/2}.
$$
Hence, by \eqref{eqcsds5} and \eqref{eqcsds4}, 
$$
\mu\bigg\{x \in M \setminus \cup_i 4B_i ; \bigg|\sum_i T_ib_i(x)\bigg| >\frac \alpha 3\bigg\} 
\le C \mu\big( \cup_i B_i \big) \le 
\frac{C}{\alpha^q}\int_M |\nabla f|^q\, d\mu.$$

It remains to handle  the term $B$.  Define $$\beta_k = \sum_{i, r_i=2^k} \frac{b_i}{r_i}
$$ for $k \in \mathbb{Z}$. With this definition, it is easy to see that 
$$
\sum_i U_ib_i = \sum_{k\in \mathbb{Z}} \int_{4^k}^\infty
 \left( \frac {2^k} {\sqrt t}\right) t\Delta e^{-t \Delta}  \beta_{k}\frac{dt}t = \int_{0}^\infty t\Delta e^{-t \Delta}  f_{t}\frac{dt}t$$
where 
$$f_{t}= \sum_{k\, ; 4^k \le t } \left( \frac {2^k} {\sqrt t}\right) \beta_{k}.
$$
By using duality from the well-known Littlewood-Paley estimate 
$$
\left\| \left( \int_{0}^\infty  |t\Delta e^{-t \Delta}  f|^2\frac{dt}t \right)^{1/2} \right\|_{q'} \le C \|f\|_{q'}
$$
(see \cite{topics}),
we find that 
$$
\left\| \sum_i U_ib_i \right\|_{q} \le C \left\| \left( \int_{0}^\infty  | f_{t}|^2\frac{dt}t\right)^{1/2}\right\|_{q}.
$$
Now,  by Cauchy-Schwarz inequality, 
$$
| f_{t}|^2 \le 2  \sum_{k\, ; 4^k \le t } \left( \frac {2^k} {\sqrt t}\right) |\beta_{k}|^2
$$
and it is easy to obtain
$$
\left\| \left( \int_{0}^\infty  | f_{t}|^2\frac{dt}t\right)^{1/2}\right\|_{q} \le C  \left\|\left(\sum_{k\in \mathbb{Z}}  |\beta_k|^2\right)^{1/2} \right\|_q.
$$
Using the bounded overlap property \eqref{eqcsds5}, one has that 
$$
\left\| \left( \sum_{k\in \mathbb{Z}}  |\beta_k|^2\right)^{1/2} \right\|_q^q \le C\int_M \sum_{i}  \frac{|b_i|^q}{r_i^q} \, d\mu
$$
and by a similar argument to one in the proof of Proposition \ref{lemmaCZD-manifold}, 
$$
\int_M \sum_{i}  \frac{|b_i|^q}{r_i^q}\, d\mu \le  C \alpha^q  \sum_{i}\mu(B_{i}).$$
Hence, by  \eqref{eqcsds4}
$$
\mu\bigg\{x \in M; \bigg|\sum_i U_ib_i(x)\bigg| 
>\frac \alpha 3\bigg\} \le  C   \sum_{i}\mu(B_{i})\le \frac{C}{\alpha^q}\int_M
|\nabla f|^q\, d\mu.$$
This concludes the proof of \eqref{eq21} when $\mu(M)=\infty$.

When $\mu(M)<\infty$, the previous argument holds for $\alpha$ such that 
$\alpha^q > \frac C {\mu(M)} \int_M |\nabla f|^q\, d\mu$. On the other hand, 
if 
$\alpha^q \le \frac C {\mu(M)} \int_M |\nabla f|^q\, d\mu$, then
\begin{equation*}
\mu\{x \in M; |\Delta^{1/2}f(x)| >\alpha\} \le  \mu(M) \le
\frac{C}{\alpha^q}\int_M |\nabla f|^q\, d\mu.
\end{equation*}

\subsection{An interpolation argument}

It is not known whether the spaces defined by the seminorms $\norm{|\nabla f|}_q$ interpolate 
by the real method. So it is not immediate to obtain $(RR_{p})$ for $q<p<2$
directly from $(E_2)$ and \eqref{weaktypeq}. We next prove this fact by adapting Marcinkiewicz theorem argument which bears again on our Calder\'on-Zygmund decomposition.

We first  do the proof when $\mu(M)=\infty$. Fix $q<p<2$ and $f\in \ccc_{0}^\infty(M)$. 
We want to show that
$$
\norm{\Delta^{1/2} f}_p\leq C_p\norm{|\nabla f|}_p. 
$$
Choose
$0<\delta<1$ so that $q<p\delta$. For $\alpha>0$, we can apply the Calder\'on-Zygmund decomposition of Proposition \ref{lemmaCZD-manifold}  with exponent $p\delta$ and threshold $\alpha$.  We may do this since $\norm{|\nabla f|}_{p\delta}<\infty$ and $(P_{p\delta})$ holds. Of course we do not want to use 
$\norm{|\nabla f|}_{p\delta}$ in a quantitative way.  We obtain that $f=g_{\alpha}+ b_{\alpha}$ with $b_{\alpha}=\sum_{i}b_{i}$.

Write 
\begin{align*}
\norm{\Delta^{1/2} f}_{p}^p &= p 2^p \int_{0}^\infty \alpha^{p-1} \mu\{ x \in M; |\Delta^{1/2}f(x)| >2\alpha\} \, d\alpha
\\
&
\le p 2^p \int_{0}^\infty \alpha^{p-1} \mu\{ x \in M; |\Delta^{1/2}g_{\alpha}(x)| >\alpha\} \, d\alpha
\\
&
\qquad + p 2^p \int_{0}^\infty \alpha^{p-1} \mu\{ x \in M; |\Delta^{1/2}b_{\alpha}(x)| >\alpha\} \, d\alpha
\\
&
\le 
I + II
\end{align*}
with
$$
I = C p 2^p \int_{0}^\infty \alpha^{p-1}  \frac {\norm{|\nabla g_{\alpha}|}_2^2} {\alpha^2}   \, d\alpha
$$
and
$$
II= C p 2^p \int_{0}^\infty \alpha^{p-1}  \frac {\norm{|\nabla b_{\alpha}|}_q^q} {\alpha^q}   \, d\alpha,
$$
where we used  $(E_2)$ and the assumption \eqref{weaktypeq}. To estimate these integrals, we need to come back to the construction of $\nabla g_{\alpha}$ and $\nabla b_{\alpha}$. 
Write $F_{\alpha}$ as the complement of $\Omega_{\alpha}= \{{\cal M}(|\nabla f|^{p\delta}) > \alpha^{p\delta }\}$. Then 
recall that $\nabla g_{\alpha}= (\nabla f) \mathbf{1}_{F_{\alpha}} + h \mathbf{1}_{\Omega_{\alpha}} $ where $|h|\le C\alpha$ and $|\nabla f |\le \alpha$ on $F_{\alpha}$. Thus $I$ splits into $I_{1}+I_{2}$ according to this decomposition. The treatment
of $I_{1}$ is done using the definition of $F_{\alpha}$, Fubini's theorem and $p<2$ as follows:
\begin{align*}
I_{1}&= \frac{C p 2^p}{2-p} \int_M    |\nabla f|^2 \ \left({\cal M}( |\nabla f|^{p\delta})\right)^{\frac {p-2}{p\delta }}   \, d\mu
\\
&
\le \frac{C p 2^p}{2-p} \int_M    |\nabla f|^p \, d\mu,
\end{align*}
where we used $|\nabla f|^2= |\nabla f|^p\ |\nabla f|^{2-p} \le |\nabla f|^p\ \left({\cal M}( |\nabla f|^{p\delta})\right)^{\frac {2-p}{p\delta }} $ almost everywhere. 
For $I_{2}$, we only use the bound of $h$ to obtain
\begin{align*}
I_{ 2}&\le  C p 2^p \int_{0}^\infty \alpha^{p-1}  \mu(\Omega_{\alpha})   \, d\alpha
\\
&
= C 2^p \int_M\left({\cal M}( |\nabla f|^{p\delta})\right)^{\frac {1}{\delta }}   \, d\mu
\\
&
\le C \int_M   |\nabla f|^p \, d\mu
\end{align*}
using the strong type $(\frac 1 \delta , \frac 1 \delta )$ of the maximal operator.

Next, we turn to the term $II$. We have $\nabla b_{\alpha}= (\nabla f) \mathbf{1}_{\Omega_{\alpha}} - h \mathbf{1}_{\Omega_{\alpha}} $ so that $II  \le 2^q(II_{1}+II_{2})$. For $II_{1}$ we have by using H\"older's inequality 
and  the strong type $(\frac 1 \delta , \frac 1 \delta )$ of the maximal operator
\begin{align*}
II_{1}&= \frac{C p 2^p}{p-q} \int_M    |\nabla f|^q \ \left({\cal M}( |\nabla f|^{p\delta})\right)^{\frac {p-q}{p\delta }}   \, d\mu
\\
&
\le \frac{C p 2^p}{p-q} \left(  \int_M    |\nabla f|^p \, d\mu\right)^{q/p} \left( \int _M   \left( {\cal M}( |\nabla f|^{p\delta})\right)^{(\frac {p-q}{p\delta })(\frac p q)'}   \, d\mu\right)^{1/(\frac p q)'}
\\
&
\le C \int _M   |\nabla f|^p \, d\mu.
\end{align*}
The treatment of the term $II_{2}$ with $h$ is as $I_{2}$. This finishes the argument when 
$\mu(M)=\infty$.

The modifications are as follows when $\mu(M)<\infty$.  The estimates apply for the part of the integral where  $\alpha>a$ with $a^{p\delta } = \frac C {\mu(M)} \int_M |\nabla f|^{p\delta} \, d\mu$. The part where $\alpha\le a$ is also bounded by $a^p \mu(M)$ which,  by H\"older inequality, is bounded  by $C \int_M    |\nabla f|^p \, d\mu$.

\section{$(R_{p})$ for  $p>2$}

In this section, we prove Theorem \ref{Rp>2} as a consequence of the next two results. 

\begin{theo}\label{newCNS}  Let $M$ be a complete non-compact Riemannian manifold satisfying $(D)$ and $(P_{2})$.  Then there exists    $p_0 \in (2,\infty]$  such that for any $q \in (2,p_{0})$ the following assertions are equivalent. 
\begin{enumerate}
  \item $(R_{p})$ holds for $2<p<q$,
 \item $(\Pi_{p})$ holds for $2<p<q$,
  \item For any $p\in (2,q)$, there exists a constant $C>0$ such that for any ball $B$ and any harmonic function  $u$ in 3B, one has the reverse H\"older inequality
$$
\left(\frac1{\mu(B)}\int_B |\nabla u|^{p}\, d\mu\right)^{\frac1{p}} \le C  \left(\frac1{\mu(2B)}\int_{2B} |\nabla u|^{2}\, d\mu\right)^{\frac1{2}} .
\eqno(RH_p)
$$
\end{enumerate}
 \end{theo}
 
 \begin{pro} \label{propRH}Let $M$ be a complete non-compact Riemannian manifold satisfying $(D)$ and $(P_{2})$.  Then there is    $p_1 \in (2,\infty]$ such that $(RH_p)$ holds for $2<p<p_{1}$.
 \end{pro}

The value of $p_{1}$ in Proposition \ref{propRH} is not known. The same is true for $p_{0}$ in Theorem \ref{newCNS}. However, if we assume $(P_{q})$ for $q\in(1,2)$ then the argument shows that $p_{0}>q'$ and for $q=1$, $p_{0}=\infty$. 

We shall first prove Proposition \ref{propRH}. Of course, harmonic functions are smooth, but the point of $(RH_p)$ is that the estimate is scale invariant.  Then we shall prove Theorem \ref{newCNS}, in establishing successively that $3. \Longrightarrow 2. \Longrightarrow 1. \Longrightarrow 3$.  This will prove $(R_{p})$ for $2<p<
\inf(p_{0}, p_{1})$.

\subsection{Reverse H\"older inequality for the gradient of harmonic functions}

Assume $(D)$ and $(P_{2})$. First we have a Caccioppoli inequality:  Let $u$ be a harmonic function on $3B$ where $B$ is some fixed ball. 
Let $B'$ be a ball such that $3B' \subset 3B$.  
Then, we have
\begin{equation}
\label{caccio}
\left(\frac1{\mu(B')}\int_{B'} |\nabla u(x)|^{2}\, d\mu\right)^{\frac1{2}} \le \frac C {r(B')} \left(\frac1{\mu(2B')}\int_{2B'} | u(x)- u_{2B'}|^{2}\, d\mu\right)^{\frac1{2}}.
\end{equation}
Its proof is entirely similar to the one in the Euclidean setting under $(D)$ and $(P_{2})$. We skip details and refer, \textit{e.g.}, to Giaquinta's book \cite{Gia}.

Next, we use Lemma \ref{selfimprovement} which tells us that $(P_{2-\ep})$ holds for some $\ep>0$. 
According to \cite{FPW}, Corollary 3.2, 
we have the $L^{2-\ep}-L^2$ Poincar\'e inequality
\begin{equation}
\label{poame}
\left(\frac1{\mu(2B')}\int_{2B'} | u(x)- u_{2B'}|^{2}\, d\mu\right)^{\frac1{2}} \le 
C r(B') \left(\frac1{\mu(2B')}\int_{2B'} | \nabla u(x) |^{2-\ep}\, d\mu\right)^{\frac1{2-\ep}}
\end{equation}
provided for any ball $B$ and subball  $B'$  
\begin{equation}
\label{FPW}
\frac {r(B')}{r(B)} \lesssim  \left( \frac {\mu(B')}{\mu(B)}\right)^{\frac 1 {2-\ep} - \frac 1 2 }\ .
\end{equation}
Admit  \eqref{FPW} and combine  \eqref{poame}  with \eqref{caccio} to obtain
a reverse H\"older inequality, 
$$
\left(\frac1{\mu(B')}\int_{B'} |\nabla u(x)|^{2}\, d\mu\right)^{\frac1{2}} \le
C \left(\frac1{\mu(2B')}\int_{2B'} | \nabla u(x) |^{2-\ep}\, d\mu\right)^{\frac1{2-\ep}}.
$$
Applying Gehring's self-improvement of reverse H\"older inequality \cite{Geh} (see also \cite{Iw}, \cite{Gia}), which holds since we work in a doubling space,  we conclude that there is
$\delta>0$ and a constant $C$ such that 
$$
\left(\frac1{\mu(B)}\int_{B} |\nabla u(x)|^{2+\delta}\, d\mu\right)^{\frac1{2+\delta}} \le
C \left(\frac1{\mu(2B)}\int_{2B} | \nabla u(x) |^{2}\, d\mu\right)^{\frac1{2}}.
$$
It remains to verify \eqref{FPW}. Write $B=B(x,r)$ and $B'=B(y,s)$ with $s<r$. Then observe that $(D)$ and $d(x,y)<r$ imply that $V(x,r)\sim V(y,r)$. Hence, we may assume that $x
=y$ and \eqref{FPW} becomes 
$$
\frac s r  \lesssim \left( \frac {V(x, s)} {V(x,r)}\right)^a$$
with $a= {\frac 1 {2-\ep} - \frac 1 2 }>0$. The doubling property $(D)$ implies that for some $\beta>0$, 
$$
 \frac {V(x, r)} {V(x,s)}  \lesssim \left (\frac r s \right)^\beta,$$ hence it suffices to have
$\beta a \le 1$. Choosing $\ep$ smaller if necessary, we obtain \eqref{FPW}. Finally, 
$(RH_{p})$ holds for $2<p<2+\delta $. 

\subsection{From  reverse H\"older  to  Hodge projection}

The main tool is  the adaptation to spaces of homogeneous type of  a result by Shen in \cite{Shen}  essentially similar to Theorem 2.1 in \cite{ACDH}. For the sake of completeness we include its proof  in  Section \ref{sec:thT}. Let $\mathcal{M}$ denote the Hardy-Littlewood maximal function.

\begin{theo} \label{thT} Let $(E,d,\mu)$ be a  measured metric space satisfying the doubling property (D).
Let $T$ be a bounded sublinear operator from $L^2(E,\mu)$ to $L^2(E,\mu)$. Assume that for  $q\in (2,\infty]$,   $1<\alpha<\beta$ and $C>0$, we have
\begin{equation}
\label{T}
\left(\frac{ 1}{\mu(B)} \int_{B} |Tf|^q\, d\mu \right)^{1/q} \le  C \left(\frac{ 1}{\mu(\alpha B)} \int_{\alpha B} |Tf|^2\, d\mu \right)^{1/2} 
\end{equation}
for all balls $B$ in $E$ and $f \in L^{2}(E, \mu)$ supported on $E\setminus \beta B$.  
Then, $T$ is bounded from  $L^p(E,\mu)$ to $L^p(E,\mu)$ for $2<p<q$. More precisely, there exists a
constant $C'$ such that for any $f\in L^p\cap L^2(E,\mu)$, we have $Tf \in L^p(E,\mu)$ and 
$$
\| Tf \|_{p} \le C' \|f\|_{p}.
$$
\end{theo}

In this statement, the functions $f$ can be vector-bundle-valued  and $|f|$ is then the norm of $f$ while $Tf$ is real valued. \bigskip

We now prove $3. \Longrightarrow 2.$ in Theorem \ref{newCNS}. We assume the reverse H\"older condition. Let  $T$ be the sublinear bounded operator  from $L^2T^*M$ into $L^2(M,\mu)$ such that 
$T\omega=|d\Delta^{-1}\delta \omega|$ when $\omega \in L^2T^*M $. Let $2<p<\tilde p< q$ where $q$ is the exponent in condition 3.
  Let $B$ be a ball in $M$ and $\omega \in L^2{T^*M}\cap L^pT^*M$ be supported on $M\setminus 4B$. Let $u$ be a distribution defined by  $\norm {|du|}_{2}<+\infty$ and  $\Delta u = \delta \omega$, so that $|du|=T\omega$. Given the support of $\omega$, it follows that  $u$ is harmonic in $3B$.  The reverse H\"older condition yields \eqref{T} with $q$ replaced by $\tilde p$, hence,  according to Theorem \ref{thT},
$$
\| T\omega \|_{p} \le C \|\omega\|_{p}.
$$
A density argument concludes the proof.

\subsection{From  Hodge projection to Riesz transform}

We begin with the proof of  Lemma \ref{reciproque}. To do this, we look at the form version of the Riesz transform,
$d\Delta^{-1/2}$, where $d$ is the exterior derivative. We assume that for $f\in \ccc_{0}^\infty(M)$
$$
\|\Delta^{1/2}f\|_{p'}\le C_{p'}\norm{|df|}_{p'}$$
and 
for $\omega\in \ccc_{0}^\infty(T^*M)$,
$$
\norm{|d\Delta^{-1}\delta \omega|}_{p} \le C_{p} \norm{\omega}_{p}. \eqno(\Pi_{p})
$$
Since $d\Delta^{-1}\delta$ is self-adjoint, the last inequality holds with 
$p$ replaced by $p'$.

Let $\omega\in \ccc_{0}^\infty(T^*M)$. Then using successively $(\Pi_{p'})$ et $(RR_{p'})$, 
$$
\norm{\Delta^{-1/2}\delta \omega}_{p'} = \norm{\Delta^{1/2}\Delta^{-1} \delta \omega}_{p'} \le C \norm{|d\Delta^{-1}\delta \omega|}_{p'} \le C \norm{\omega}_{p'}.$$
Hence, by duality, $d\Delta^{-1/2}$ is bounded on $L^p$.

\bigskip

The proof that $2. \Longrightarrow 1.$ in Theorem \ref{newCNS} is now easy. By combining Theorem \ref{poincareq} with Lemma \ref{selfimprovement}, we have $(RR_{p})$ for $2-\ep<p<2$.  Let $p_{0}=(2-\ep)'$ and $2<q<p_{0}$. If we assume $(\Pi_{p})$  for $2<p<q$, then Lemma \ref{reciproque} gives us $(R_{p})$ for $2<p<q$. 

\subsection{From Riesz transform  to  reverse H\"older inequalities}

We show here the necessity of the reverse H\"older inequalities $ (RH_{p})$.
We assume that the Riesz transform is bounded on $L^p$ for $2<p<q$. Fix such a $p$.

Let  $B$ be a ball, $r$ its radius and let $u$ be harmonic function  in $3B$. Let $\varphi$ a $C^1$ function, supported in $2B$ with $\varphi=1$ on $\frac 3 2 B$, $\|\varphi\|_{\infty}\le 1$ and $\|\nabla \varphi\|_{\infty} \le C/r$. We assume that $\int_{2B}u=0$ so that  it follows from $(P_{2})$ that 
$$
 r^{-2}\int_{2B} |u|^2 d\mu + \int_{2B}|\nabla (u\varphi)|^2\, d\mu \le C \int_{2B}|\nabla u|^2\, d\mu.
 $$
To estimate $\int_{B} |\nabla u |^p d\mu$, it suffices to estimate $\int_{B} |\nabla (u\varphi)|^p \, d\mu$. Using   an idea in  \cite{AT}, p.~35, we can write 
$$
u\varphi= e^{-r^2\Delta}(u\varphi) + u\varphi - e^{-r^2\Delta}(u\varphi) = e^{-r^2\Delta}(u\varphi) - \int_{0}^{r^2} e^{-s\Delta}\Delta(u\varphi) \, ds,
$$
hence
$$
\nabla (u\varphi)= \nabla e^{-r^2\Delta}(u\varphi) - \int_{0}^{r^2} \nabla e^{-s\Delta}\Delta(u\varphi) \, ds.
$$

Let $p<\rho<q$. Since the Riesz transform is bounded on $L^\rho$, by the easy part of the necessary and sufficient condition in Theorem \ref{nsc}, we have that $\sqrt t \nabla e^{-t\Delta}$ is bounded on $L^\rho$ uniformly with respect to $t$.  It essentially follows from  Lemma 3.2 in \cite{ACDH}   that \begin{equation}\label{eqoffdiag}
 \left( \frac 1{\mu(B)}\int_B   |\nabla e^{-s\Delta}f|^p\, d\mu \right)^{1/p}   \le \frac {Ce^{-\frac{\alpha 4^j r^2}{s}}} {\sqrt s} \left(\frac 1 {\mu(c_{2}2^{j}B)} \int_{C_{j}(B)} |f|^2\, d\mu\right)^{1/2}
\end{equation}
for some constants $C$ and $\alpha$ depending only on $(D)$, $(P_{2})$,  $p$ and $\rho$ whenever $f$ is supported in $C_{j}(B)$ and $s \lesssim r^2(B)$. Here $C_{1}(B)$ is a fixed multiple of $B$, and for $j\ge 2$, 
$C_{j}(B)$ is a ring based on $B$: there are constants $c_{1},c_{2}$ such that for all $j\ge 1$,  if $x\in C_{j}(B)$ then $c_{1}2^jr \le d(x,B)  \le c_{2}2^jr$.  

It suffices to apply this inequality to $f=u\varphi$ which is supported in $2B$ to treat  the $L^p$ average of $\nabla e^{-r^2\Delta}(u\varphi) $ on $B$.
 
In the other term, a computation yields $$\Delta(u\varphi)= - d u\cdot d \varphi - \delta ( u  d \varphi ).$$ We replace $\Delta(u\varphi)$ by its expression and observe that the support
condition of $d\varphi$ allows us to use  the previous estimates \eqref{eqoffdiag}  for 
$\nabla e^{-s\Delta}(du\cdot d\varphi)$ when $j\ge 2$.  Then, by Minkowski inequality, 
\begin{equation*}
\left(\frac1{\mu(B)}\int_B \left| \int_{0}^{r^2} \nabla e^{-s\Delta}(du\cdot d\varphi) \, ds\right |^{p}d\mu\right)^{\frac1{p}} \le   C \left(\frac 1{\mu(2B)}\int_{2B} |\nabla u|^{2}\, d\mu\right)^{\frac1{2}} .
\end{equation*}

For the remaining term, it suffices to prove 
\begin{equation}\label{eqoffdiag2}
 \left( \frac 1{\mu(B)}\int_B   |\nabla e^{-s\Delta}\delta f|^p\, d\mu \right)^{1/p}   \le \frac {Ce^{-\frac {c r^2}s}}{s} \left(\frac 1 {\mu(2B)} \int_{2B \setminus \frac 3 2 B} |f|^2\, d\mu\right)^{1/2}
\end{equation} whenever $f$ is supported in $2B \setminus \frac 3 2 B$ and $s\le r^2$ since this yields
\begin{equation*}
\left(\frac1{\mu(B)}\int_B \left| \int_{0}^{r^2}\nabla e^{-s\Delta}\delta(ud\varphi) \, ds\right |^{p}d\mu\right)^{\frac1{p}}  \le  \frac C r \left(\frac 1{\mu(2B)}\int_{2B} | u |^{2}\, d\mu\right)^{\frac1{2}},
\end{equation*}
which concludes the proof of  $(RH_p)$.

To see  \eqref{eqoffdiag2}, the strategy is as follows. We use that $\nabla e^{-t\Delta}\delta 
= (\nabla e^{-t/2\Delta}) (e^{-t/2\Delta} \delta) $. For the second operator we have the Gaffney type 
estimate 
\begin{equation*}\label{Gaffney}
\| \sqrt t \, e^{-t\Delta}\delta \omega \|_{L^2(F)} 
 \le C e^{-\frac {\alpha d(E,F)^2} t} \|\omega\|_{L^2(E)}.
\end{equation*}
whenever $f$ is a 1-form supported on $E$ and $E, F$ are closed subsets of $M$ and $t>0$. This estimate is for example proved in \cite{ACDH} for the dual operator $de^{-t\Delta}$.
Make use  of it with $E= 2B \setminus \frac 3 2 B$ and successively
$F= \frac 5 4 B$, $4B \setminus  \frac 5 4 B$, and $2^{j+1}B \setminus 2^jB$ for $j\ge 2$ and combine them with \eqref{eqoffdiag} to conclude. Similar calculations
are shown in \cite{ACDH} and we skip further details.

\section{Proof of Theorem \ref{thT}} \label{sec:thT} We split the argument in several steps.
The following  lemma is a localisation result and  is applied in the proof of a good lambda inequality which is the key step. The latter yields $L^p$ inequalities, which applied to our particular hypotheses concludes the proof.

\begin{lem} (localisation) There is $K_0$ depending only on the doubling constant of
$E$ such that the following holds. Given 
$f\in L^1_{loc}(E,\mu)$,  a ball
$B$ and  $\lambda>0$ such that there exists $\bar x \in B$ for which
$\mathcal{M}f(\bar x)
\le
\lambda$, then for any $K\ge K_0$,
$$
\{\chi_B\mathcal{M}f >K\lambda\} \subset  \{\mathcal{M}(f\chi_{3B}) >\frac K {K_0}\lambda\}.
$$ 
\end{lem}

\paragraph{Proof:} Recall that $\mathcal{M}$ is comparable to the centered maximal function
$\mathcal{M}_c$: there is $K_0$ depending only on the doubling constant such that
$\mathcal{M} \le K_0 \mathcal{M}_c$. 

Let $x \in B$ with $\mathcal{M}f(x)>K\lambda$. Then  $\mathcal{M}_cf(x) > \frac K {K_0}
\lambda$. Hence, there is a ball $B(x,r)$ centered at $x$ with radius $r$ such that
$$
\frac {1}{\mu(B(x,r))}\int_{B(x,r)} f\, d\mu > \frac K {K_0}
\lambda.
$$
If $\frac K {K_0}
\ge1$,  $\bar x \notin B(x,r)$ since $\mathcal{M}f(\bar x) \le \lambda$. 
The conditions $x\in B$, $\bar x \in B$ and $\bar x \notin B(x,r)$ imply
$B(x,r) \subset 3B$. Hence,
$$
 \frac K {K_0}
\lambda < \frac {1}{\mu(B(x,r))}\int_{B(x,r)} (f\chi_{3B})\, d\mu \le
\mathcal{M}(f\chi_{3B})(x).
$$  

This proves the lemma. 

\bigskip

We continue with a two parameters family of good lambda inequalities.

\begin{pro} (two parameter good-lambda inequalities) Fix $1<q\le \infty$ and $a>1$.  Let $F,G\in L^1_{loc}(E,\mu)$, non-negative. We say that  $(F,G) \in \mathcal{E}_{q,a}$ if  one can find for every
ball $B$ non-negative measurable functions $G_B,H_B $ defined on $B$ with 
$$
F\le G_B+H_B \quad a.e. \ on \ B
$$ such that 
$$
\left( \frac {1}{\mu(B)}\int_{B} (H_B)^{q}\,
d\mu\right)^{1/q} \le a \inf _{x\in B} \mathcal{M}F(x) + \inf_{x\in B}, G(x)
$$
$$
 \frac {1}{\mu(B)}\int_{B} G_B\,
d\mu  \le \inf_{x\in B} G(x).
$$
There exist $C=C(q,
(D), a)$ and $K'_0=K'_0(a,(D))$ such that for  $(F,G) \in \mathcal{E}_{q,a}$,
 for all $\lambda>0$, for all $K>K_0'$ and $\gamma\le 1$, 
$$
\mu\{\mathcal{M}F>K\lambda, G\le \gamma\lambda\} \le C\left(\frac{1}{K^q} + \frac
\gamma K\right) \mu\{\mathcal{M}F>\lambda\}
$$
provided $\{\mathcal{M}F>\lambda\}$ is a proper subset of $E$.

If $q=\infty$, we understand the average in $L^q$ as an essential supremum. In this
case, we set $\frac{1}{K^q} =0$.
\end{pro}

\paragraph{Proof:} Let $E_\lambda=\{\mathcal{M}F>\lambda\}$. This is an open proper subset of $E$.  The Whitney  decomposition 
for $E_\lambda$
yields a family of boundedly overlapping balls $B_i$ 
such that $E_\lambda=\cup_iB_i$. There exists
$c>1$ such that, for all
$i$, $cB_i$ contains at least one point $\overline{x_i}$ outside
$E_\lambda$, that is
$${ \mathcal{M}}F(\overline{x_i})\le \lambda.$$

Let 
$B_\lambda=\{\mathcal{M}F>K\lambda, G\le \gamma\lambda\}$. If $K\ge 1$ then 
$B_\lambda \subset E_\lambda$, hence
$$
\mu(B_\lambda) \le \sum_i \mu(B_\lambda \cap B_i) \le \sum_i \mu(B_\lambda
\cap cB_i).
$$
Fix $i$. If $B_\lambda
\cap cB_i=\emptyset$, we have nothing to do. If not, there is a point
$\overline{y_i} \in cB_i$ such that 
$$G(\overline{y_i}) \le
\gamma\lambda.$$
By the localisation lemma applied to $F$ on $cB_i$, if $K\ge K_0$, then
$$
\mu(B_\lambda
\cap cB_i) \le \mu (\{\mathcal{M}F>K\lambda\} \cap cB_i) \le
\mu\{\mathcal{M}(F\chi_{3cB_i})>
\frac K {K_0} \lambda\}.
$$
Now use $F\le G_i+H_i$ on $3cB_{i}$ with $G_i=G_{3cB_i}$ and $H_i=H_{3cB_i}$ to deduce
$$\mu\{\mathcal{M}(F\chi_{3cB_i})>
\frac K {K_0} \lambda\} \le 
\mu\{\mathcal{M}(G_i\chi_{3cB_i})>
\frac K {2K_0} \lambda\} + \mu\{\mathcal{M}(H_i\chi_{3cB_i})>
\frac K {2K_0} \lambda\}.
$$
Now by using the weak type $(1,1)$ and $(q,q)$ of the maximal operator
with respective constants $c_1$ and $c_q$, we have
$$
\mu\{\mathcal{M}(G_i\chi_{3B_i})>
\frac K {2K_0} \lambda\} \le \frac {2K_0c_1}{K\lambda} \int_{3cB_i} G_i\,
d\mu \le \frac {2K_0c_1}{K\lambda} \mu(3cB_i) G(\overline{y_i}) \le 
\frac {2K_0c_1 \gamma}{K} \mu(3cB_i),
$$
and, if $q<\infty$,
\begin{align*}
\mu\{\mathcal{M}(H_i\chi_{3cB_i})>
\frac K {2K_0} \lambda\} &\le \left(\frac {2K_0c_q}{K\lambda}\right)^q
\int_{3cB_i} H_i^q\, d\mu\\
& \le \left(\frac {2K_0c_q}{K\lambda}\right)^q
\mu(3cB_i) (a \mathcal{M}F(\overline{x_i}) + G(\overline{y_{i}})) ^q 
\\
&\le 
\left(\frac {2K_0c_q(a+1)}{K}\right)^q \mu(3cB_i) .
\end{align*}
Hence, summing over $i$ yields
$$\mu(B_\lambda) 
\le 
C\left(\frac{1}{K^q} + \frac
\gamma K\right) \sum_{i}\mu(3cB_i) \le C'\left(\frac{1}{K^q} + \frac
\gamma K\right) \mu(E_\lambda) 
$$
by applying the doubling property together with the bounded overlap.
If $q=\infty$, then 
$$
\|\mathcal{M}(H_i\chi_{3cB_i})\|_\infty \le  \|H_i\chi_{3cB_i}\|_\infty  \le a
\mathcal{M}F(\overline{x_i}) + G(\overline{y_{i}})\le (a +1)\lambda,
$$
so that, choosing $K\ge 2K_0(a+1)$ leads us to $\{\mathcal{M}(H_i\chi_{3B_i})>
\frac K {2K_0} \lambda\}=\emptyset$. The rest of the proof is unchanged. 
This proves the proposition.

\begin{cor} Assume that $(F,G) \in \mathcal{E}_{q,a}$.  Let $1<\rho<q$ and assume that $\|G\|_{\rho}<\infty$ and    $\|F\|_{1}<\infty$. Then, we have 
$$
\|\mathcal{M}F\|_\rho \le C\left( \|G\|_\rho +\mu(E)^{\frac 1 \rho - 1}  \|F\|_{1}\right),\footnote{In the case $\mu(E)=\infty$, the last term vanishes but we still need some a priori knowledge such as 
$F\in L^1$ to conclude.}
$$
where the constant $C$ depends on $(D)$, $\rho$, $q$, $a$.
\end{cor}

\paragraph{Proof:} We begin with the case $\mu(E)=\infty$. Define $\Phi(t)= p\int_{0}^t \lambda^{\rho-1}\mu\{\mathcal{M}F>\lambda\}\, d\lambda$ for $t\ge 0$. Since  $\|F\|_{1}<\infty$, the maximal theorem implies that   $\lambda \mu\{\mathcal{M}F>\lambda\}$ is bounded  on $\mathbb{R}^+$. As $1<\rho$,  $\Phi$ is a well-defined positive and non-decreasing function on $\mathbb{R}^+$ into $\mathbb{R}^+$. 

By the maximal theorem and $\|F\|_{1}<\infty$, $\{\mathcal{M}F>\lambda\}$ is a proper subset  in $E$,  hence the good lambda inequality is valid and integration  leads us to 
$$
\Phi(Kt) \le CK^\rho\left(\frac{1}{K^q} + \frac
\gamma K\right) \Phi(t) + \left(\frac K \gamma\right)^\rho \|G\|_\rho^\rho.
$$
Since $\rho<q$, one can choose $K$
large enough and $\gamma$ small enough so that 
$$
CK^\rho\left(\frac{1}{K^q} + \frac
\gamma K\right) \le \frac 1 2.
$$
hence, for this choice, for all $t\ge 0$
$$
\Phi(Kt) \le \frac 1 2 \Phi(t) +  \left(\frac K \gamma\right)^\rho \|G\|_\rho^\rho.
$$
An easy iteration proves that $\Phi$ is bounded and this  proves the corollary in this case
as $\Phi(\infty)$ is  $\|\mathcal{M}F\|_{\rho}^\rho$.

In the case where  $\mu(E)<\infty$, we have $\lambda \mu\{\mathcal{M}F>\lambda\} \le C \|F\|_{1}$, hence for $\lambda>a$ with $a= \frac C {\mu(E)} \|F\|_{1}$, the good lambda
inequality applies. If we define $\Phi$ as before, the previous argument gives us a control of 
$\Phi(\infty)- \Phi(a)$ by $C\|G\|_{\rho}^\rho$ and it remains to controlling $\Phi(a)$.  But 
$\Phi(a) \le a^\rho \mu(E)$ and the conclusion follows. 
\bigskip

Now, we may prove Theorem \ref{thT}. We let $f \in L^p\cap L^2(E,\mu)$ and   $F=|Tf|^2$. We let $G_{B}=2|T(\chi_{\beta B}f)|^2$ and $H_{B}= 2|T( (1-\chi_{\beta B})f)|^2$. On the one hand,  for $C$  depending only on $(D)$ and the norm $\|T\|$ of $T$ on $L^2$,
$$
 \frac {1}{\mu(B)}\int_{B} G_B\,
d\mu  \le  \frac {2\|T\|^2}{\mu(B)} \int_{\beta B} |f|^2 \le  C \inf_{x\in B}  \mathcal{M}(|f|^2)(x).
$$
On the other hand, since $  (1-\chi_{\beta B})f$ is supported away from $\beta B$, the assumption \eqref{T} yields
$$
\left( \frac {1}{\mu(B)}\int_{B} (H_B)^{q/2}\,
d\mu\right)^{2/q} \le  \frac {C}{\mu(\alpha B)}\int_{\alpha B} H_B\,
d\mu 
$$
and we have 
$$
\int_{\alpha B} H_B\,
d\mu  \le 4\int_{\alpha B} F\,
d\mu  + 2\int_{\alpha B} G_B\,
d\mu  
$$
hence for some $a>0$, 
$$
\left( \frac {1}{\mu(B)}\int_{B} (H_B)^{q/2}\,
d\mu\right)^{2/q}   \le a \inf _{x\in B} \mathcal{M}F(x) + C \inf_{x\in B}  \mathcal{M}(|f|^2)(x).
$$
Thus we conclude with $G= C \mathcal{M}(|f|^2)$  that 
if $2<p<q$,  since $Tf \in L^2$ hence $F \in L^{1}$, then
$$\|F\|_{p/2 } \le C\left(\| G\|_{p/2} + \mu(E)^{\frac 2 p - 1}  \| F\|_{1}\right).$$
Observe then that $\| G\|_{p/2} \sim \|f\|_{p}^2$ and by the $L^2$ boundedness of $T$ and H\"older inequality, 
 $$\mu(E)^{\frac 2 p - 1} \| F\|_{1}  \le  C  \mu(E)^{\frac 2 p - 1} \|f\|_{2}^2 \le C \|f\|_{p}^2.$$
 
 \bigskip

{\bf Acknowledgement:} This work was triggered by a question of Theo Sturm to the second author
after a talk he gave in Banff in  April 2004. The two authors would like to thank him for this.

\bigskip

\end{document}